\documentclass{amsart}

\usepackage{amssymb}
\usepackage[english,francais]{babel}

\newtheorem{theorem}{Theorem}[section]

\newtheorem{e-proposition}[theorem]{Proposition}

\newtheorem{e-definition}[theorem]{Definition\rm}

\newtheorem{theoreme}{Th\'eor\`eme}[section]

\def\og{\leavevmode\raise.3ex\hbox{$\scriptscriptstyle\langle\!\langle$~}}
\def\fg{\leavevmode\raise.3ex\hbox{~$\!\scriptscriptstyle\,\rangle\!\rangle$}}

\usepackage[utf8]{inputenc}
\usepackage{amsmath,mathrsfs}

\newcommand{\DD}{{\mathbb D}}
\newcommand{\FF}{{\mathbb F}}
\newcommand{\NN}{{\mathbb N}}
\newcommand{\PP}{{\mathbb P}}
\newcommand{\QQ}{{\mathbb Q}}
\newcommand{\RR}{{\mathbb R}}
\newcommand{\XX}{{\mathbb X}}
\newcommand{\ZZ}{{\mathbb Z}}
\newcommand{\G}{{\mathscr G}}
\renewcommand{\H}{{\mathscr H}}
\newcommand{\OOO}{{\mathscr O}}
\newcommand{\ppp}{{\mathfrak p}}
\newcommand{\Aut}{\operatorname{Aut}}
\newcommand{\bigO}{\operatorname{O}}
\newcommand{\bs}{\backslash}
\newcommand{\card}{{\operatorname{Card}}}
\newcommand{\dbs}{\backslash\!\!\backslash}
\newcommand{\ga}{\gamma}
\newcommand{\Ga}{\Gamma}
\newcommand{\gengeod}{\operatorname{\widecheck{\G\,}\!\!}}
\newcommand{\GL}{\operatorname{GL}}
\newcommand{\haar}{\operatorname{Haar}}
\newcommand{\mmm}{{\mathfrak m}}
\newcommand\normalout{\partial^1_{+}}

\newcommand\normalmp{\partial^1_{\mp}}
\newcommand{\PGL}{\operatorname{PGL}}
\newcommand{\ra}{\rightarrow}
\newcommand{\redn}{\operatorname{\tt N}}
\newcommand{\vol}{\operatorname{vol}}
\newcommand{\Vol}{\operatorname{Vol}}
\newcommand{\weakstar}{\overset{*}\rightharpoonup}
\newcommand{\wt}[1]{{\widetilde{#1}}}

\usepackage{mathtools}
\makeatletter
\DeclareRobustCommand\widecheck[1]{{\mathpalette\@widecheck{#1}}}
\def\@widecheck#1#2{%
    \setbox\z@\hbox{\m@th$#1#2$}%
    \setbox\tw@\hbox{\m@th$#1%
       \widehat{%
          \vrule\@width\z@\@height\ht\z@
          \vrule\@height\z@\@width\wd\z@}$}%
    \dp\tw@-\ht\z@
    \@tempdima\ht\z@ \advance\@tempdima2\ht\tw@ \divide\@tempdima\thr@@
    \setbox\tw@\hbox{%
       \raise\@tempdima\hbox{\scalebox{1}[-1]{\lower\@tempdima\box
\tw@}}}%
    {\ooalign{\box\tw@ \cr \box\z@}}}
\makeatother

\title{Equidistribution non-archimédienne \\
et actions de groupes sur les arbres} 

\author{Anne Broise-Alamichel \and Jouni Parkkonen \and Fr\'ed\'eric Paulin}

\begin{document}

\maketitle

\selectlanguage{francais}

\medskip
\selectlanguage{francais}

\begin{abstract}
  \selectlanguage{francais} Nous donnons des résultats
  d'équidistribution d'éléments de corps de fonctions sur des corps
  finis, et d'irrationnels quadratiques sur ces corps, dans leurs corps
  locaux complétés. Nous déduisons ces résultats de théorèmes
  d'équidistribution de perpendiculaires communes dans des quotients
  d'arbres par des réseaux de leur groupe d'automorphismes, démontrés
  à l'aide de propriétés ergodiques du flot géodésique discret.
%{\it Pour citer cet article~: A. Nom1, A. Nom2, C. R.
%Acad. Sci. Paris, Ser. I 340 (2005).}
\vskip 0.5\baselineskip

\selectlanguage{english}

\noindent
\noindent{\sc Abstract. }{\bf Non-Archimedean equidistribution and group actions on trees. }
We give equidistribution results of elements of function fields over
finite fields, and of quadratic irrationals over these fields, in
their completed local fields. We deduce these results from
equidistribution theorems of common perpendiculars in quotients of
trees by lattices in their automorphism groups, proved by using
ergodic properties of the discrete geodesic flow.
%{\it To cite this article: A. Nom1, A. Nom2, C. R.
%Acad. Sci. Paris, Ser. I 340 (2005).}
\end{abstract}

\selectlanguage{english}
\section*{Abridged English version}

Let $\FF_q$ be a finite field with $q$ elements. Let $K$ be the
function field of a geometrically irreducible smooth projective curve
${\bf C}$ over $\FF_q$ of genus $g$ and let $v$ be a (normalised
discrete) valuation of $K$, with associated absolute value $|\cdot|_v
={q_v}^{-v(\cdot)}$. Let $K_v$ be the completion of $K$ for $v$, with
residual field of order $q_v$. Let $R_v$ be the affine function ring
associated to $v$. Let $\zeta_K$ be Dedekind's zeta function of $K$
and $\haar_{K_v}$ the usual Haar measure of $K_v$. Let $\Delta_x$ be
the unit Dirac mass at any point $x$ of a topological space.

%For instance, if $K=\FF_q(Y)$ is the field of rational fractions in
%one variable $Y$ over $\FF_q$ and $v=v_\infty:\frac{P}{Q}\mapsto \deg
%Q-\deg P$, then $g=0$, $q_v=q$, $K_v=\FF_q((Y^{-1}))$ is the field of
%formal Laurent series in $Y^{-1}$ over $\FF_q$,
%$\OOO_v=\FF_q[[Y^{-1}]]$ is the ring of formal Laurent series in
%$Y^{-1}$ over $\FF_q$ and $R_v=\FF_q[Y]$ is the ring of polynomials in
%$Y$ over $\FF_q$.

The aim of this note is to give equidistribution results in $K_v$ of
elements of $K$, and of quadratic irrationals in $K_v$ over $K$,
when considered in an orbit of the modular group $\PGL_2(R_v)$ acting
by homographies on $\PP_1(K_v)=K_v\cup\{\infty\}$. We only give in the
abridged English version two of these results, and we refer to
\cite{BroParPau16} for more general versions and complete proofs.

The first one, analogous to a result of Mertens on the
equidistribution of $\QQ$ in $\RR$, says that every orbit of a finite
index subgroup (not necessarily a congruence one) of the modular group
$\PGL_2(R_v)$ of an element of $K$ equidistributes, when the absolute
value of the denominators in reduced form converge to infinity.

\medskip
\begin{theorem}\label{theo:mertensEN}  
  For every finite index subgroup $G$ of
  $\GL_2(R_v)$, as $s\to+\infty$, with $G_{(1,0)}$ the stabiliser of
  $(1,0)$ in $G$, we have
$$
\frac{(q_v^2-1)\;(q_v+1)\;\zeta_K(-1)\; [\GL_2(R_v):G]}
{q_v^3\;q^{\,g-1}\;[\GL_2(R_v)_{(1,0)}:G_{(1,0)}]}
\;s^{-2} \sum_{(x,y)\in G (1,0),\; |y|_v\leq s}
\Delta_{\frac xy}\;\;\weakstar\;\; \haar_{K_v}\;.
$$
\end{theorem}

Assume now that the characteristic of $K$ is not $2$. If $\alpha\in
K_v$ is a quadratic irrational over $K$, let $\alpha^\sigma$
be its Galois conjugate and
$$
h(\alpha)=\frac{1}{|\alpha-\alpha^\sigma|_v}\;,
$$
which is an appropriate complexity in a given orbit of the modular
group $\PGL_2(R_v)$.

\medskip
\begin{theorem}\label{theo:quadratirratEN}
  For every finite index subgroup $G$ of
  $\GL_2(R_v)$ and every quadratic irrational $\alpha_0\in K_v$ over
  $K$, as $s\to+\infty$, we have
$$ \frac{(q_v+1)^2\;\zeta_K(-1)\; m_0\;[\GL_2(R_v):G]}
  {2\;q_v^2\;(q-1)\;|v(\operatorname{tr} g_0)|} \;s^{-1}
  \sum_{\alpha\in G\cdot \alpha_0,\; h(\alpha)\leq s}
  \Delta_{\alpha}\;\;\weakstar\;\; \haar_{K_v}\;,
$$
where $g_0\in G$ fixes $\alpha_0$ with $v(\operatorname{tr} g_0)\neq 0$ and
$m_0$ is the index of ${g_0}^\ZZ$ in the stabiliser of $\alpha_0$ in $G$.
\end{theorem}

\medskip Our approach relies on techniques of ergodic theory for the
geodesic flow on trees. Let $\XX$ be a $(q+1)$-regular tree, let $\Ga$
be a lattice in $\Aut(\XX)$, preserving the set of vertices at even
distance from a given vertex. Let $\Vol(\Ga\dbs\XX)=\sum_{[x]\in
  \Ga\bs V\XX} \frac{1}{|\Ga_x|}$. Let $\gengeod \XX$ be the
Bartels-Lück space of continuous maps $\ell$ from $\RR$ to the
geometric realisation $X$ of $\XX$ such that $\ell(0)$ is a vertex,
that are isometric on a closed interval, and constant on each
complementary component, endowed with the geodesic flow $(t,\ell)
\mapsto \{s\mapsto\ell(s+t)\}$ (with discrete time $t\in\ZZ$).  Let
$\DD^{\pm}$ be two subtrees of $\XX$ such that the families
$(\ga\DD^\pm)_{\ga\in\Ga/\Ga_{\DD^\pm}}$ are locally finite in $\XX$.
The closed subspace of $\gengeod \XX$ consisting of geodesic rays
arriving at $\DD^+$ (respectively exiting from $\DD^-$) carries a
natural Borel measure $\wt\sigma^-_{\DD^+}$ (respectively
$\wt\sigma^+_{\DD^-}$).

Our main result is a simultaneous equidistribution result of common
perpendiculars: for every $\ga\in\Ga$, let
$\alpha^-_\ga:[0,d(\DD^-,\ga\DD^+)]\ra X$, and
$\alpha^+_\ga:[-d(\DD^-,\ga\DD^+),0]\ra X$ be the parametrisations with
$\alpha^-_\ga(0)\in \DD^-$ and $\alpha^+_\ga(0)\in \DD^+$ of the common
perpendicular from $\DD^-$ to $\ga\DD^+$, when it exists.

\medskip
\begin{theorem}\label{theo:arbreEN} 
  As $t\ra+\infty$, for the weak-star convergence of measures
  on $\gengeod \XX\times\gengeod \XX$, we have
$$
\frac{(q^2-1)\;(q+1)}{2\;q^2}\;\Vol(\Ga\dbs\XX)\;q^{-t}
\sum_{\substack{\ga\in\Ga/\Ga_{\DD^+}\\0<d(\DD^-,\ga\DD^+)\leq t}}\;
\Delta_{\alpha^-_\ga}\otimes \Delta_{\ga^{-1}\alpha^-_\ga}
\;\;\weakstar\;\;\wt\sigma^+_{\DD^-}\otimes \wt\sigma^-_{\DD^+}\;.
$$
\end{theorem}

\medskip When $\Ga$ is geometrically finite, there is an error term in
$\bigO(q^{-\kappa})$ for some $\kappa>0$ in this equidistribution
claim evaluated on a locally constant function with compact support.
The proof (see \cite{BroParPau16}) uses the mixing property of the
square of the geodesic flow, and its exponential mixing property
announced in \cite{Kwon15} when $\Ga$ is geometrically finite.

The proof of the above arithmetic applications uses for $\XX$ the
Bruhat-Tits building $\XX_v$ of $(\PGL_2,K_v)$, on which the modular
group $\PGL_2(R_v)$ is a geometrically finite lattice. For $\DD^-$ and
$\DD^+$ in Theorem \ref{theo:quadratirratEN}, we take the same
horoball $\H_\infty$ centered at the point $\infty$ of the space of
ends $\partial_\infty\XX_v=\PP_1(K_v)$.  For $\DD^-$ and $\DD^+$ in
Theorem \ref{theo:quadratirratEN}, we take $\DD^-=\H_\infty$ and 
$\DD^+$ the geodesic line in $\XX_v$ with points at infinity
$\alpha_0$ and $\alpha_0^\sigma$.

\selectlanguage{francais}
\section{\'Equidistribution dans des corps locaux 
non-archimédiens}
\label{sec:corplocal}

Une motivation pour notre travail est le résultat de Mertens suivant,
précisant quantitativement la densité du corps $K=\QQ$ des fractions
de l'anneau $\ZZ$ dans sa complétion $\RR$ pour la valeur absolue
usuelle.  Nous noterons $\Delta_x$ la masse de Dirac unité en tout point
$x$ d'un espace topologique, $H_x$ le stabilisateur de tout point $x$
d'un ensemble muni d'une action d'un groupe $H$ et $\;\weakstar\;$ la
convergence vague des mesures sur tout espace localement compact. En
notant $\haar_\RR$ la mesure de Lebesgue de $\RR$, quand $s\ra
+\infty$, nous avons (voir par exemple \cite{ParPau14AFST} pour une
démonstration géométrique)
$$
\frac{\pi^2}{6}\;s^{-2}\sum_{p,q\in\ZZ,\;(p,q)=1,\;|q|\leq s}
\Delta_{\frac{p}{q}}\;\;\weakstar\;\;\haar_\RR\;.
$$
Le but de cette note est de donner un analogue à ce résultat
d'équidistribution de rationnels dans le cas des corps de fonctions,
ainsi que des résultats d'équidistribution d'irrationnels
quadratiques. Nous renvoyons à \cite{BroParPau16} pour des énoncés et
démonstrations complets.

\medskip Soit $\FF_q$ un corps fini d'ordre $q$. Soient $K$ le corps
des fonctions d'une courbe projective lisse géométri\-que\-ment
irréductible ${\bf C}$ sur $\FF_q$ de genre $g$, et $v$ une valuation
(discrète, normalisée) sur $K$. Soit $K_v$ la complétion de $K$ à la
place $v$, d'anneau de valuation $\OOO_v$, de corps résiduel d'ordre
$q_v$, et de valeur absolue $|\cdot|_v ={q_v}^{-v(\cdot)}$.  Soit
$R_v$ l'anneau des fonctions de la courbe affine ${\bf C} -
\{v\}$. Notons $\zeta_K$ la fonction zéta de Dedekind de $K$ et
$\haar_{K_v}$ la mesure de Haar sur le groupe additif $K_v$ normalisée
pour que $\haar_{K_v}(\OOO_v)=1$.

Notre premier résultat est un résultat d'équidistribution, analogue à
celui de Mertens, de l'orbite du point à l'infini $\infty=[1:0]$ par
un sous-groupe d'indice fini (pas forcément de congruence) du groupe
modulaire $\PGL_2(R_v)$ agissant sur les points rationnels sur $K$ de
$\PP_1(K_v)=K_v\cup\{\infty\}$.

\medskip
\begin{theoreme}\label{theo:mertens}  
  Pour tout sous-groupe d'indice fini $G$ de
  $\GL_2(R_v)$, quand $s\to+\infty$, nous avons
$$
\frac{(q_v^2-1)\;(q_v+1)\;\zeta_K(-1)\; [\GL_2(R_v):G]}
{q_v^3\;q^{\,g-1}\;[\GL_2(R_v)_{(1,0)}:G_{(1,0)}]}
\;s^{-2} \sum_{(x,y)\in G (1,0),\; |y|_v\leq s}
\Delta_{\frac xy}\;\;\weakstar\;\; \haar_{K_v}\;.
$$
\end{theoreme}

\medskip Nous renvoyons à \cite{BroParPau16} pour l'énoncé général
d'équidistribution dans $\PP_1(K_v)$ de l'orbite de tout point de
$\PP^1(K)$ par $G$ (intéressant lorsque le nombre de classes de $K$
n'est pas $1$), dont le résultat suivant se déduit: si $\mmm$ est un
idéal fractionnaire non nul de $R_v$, de norme $\redn(\mmm)$, il
existe $r_\mmm\in\{1,\dots, q-1\}$ explicite et $\kappa>0$ tels que,
quand $s\to+\infty$, pour l'action par transvections $k\cdot (x,y)=
(x+ky,y)$ de $R_v$ sur $\mmm\times\mmm$,
\begin{align*}
& \card\;\;R_v\backslash
\big\{(x,y)\in\mmm\times\mmm\;:\;0<\frac{\redn(y)}{\redn(\mmm)} 
\leq s, \;\langle x,y\rangle=\mmm\big\}\\ =\; &
\;\frac{(q-1)\;q^{\,2g-2}\;q_v^3}{(q_v^2-1)\;(q_v+1)\;\zeta_K(-1)\;r_\mmm}
\;s^{2} +\bigO(s^{2-\kappa})\;.
\end{align*}

Donnons maintenant un échantillon de nos résultats d'équidistribution
d'irrationnels quadratiques (voir \cite {BroParPau16} pour des énoncés
et démonstrations complets), en supposant la caractéristique
différente de $2$. Si $\alpha\in K_v$ est un irrationnel quadratique
sur $K$, notons $\alpha^\sigma$ son conjugué de Galois, ${\tt n}
(x-y\alpha)= (x-y\alpha)(x-y\alpha^\sigma)$ pour tous les $x,y\in K$
la forme norme associée, et
$$
h(\alpha)=\frac{1}{|\alpha-\alpha^\sigma|_v}\;,
$$
qui, comme nous allons le voir, est une complexité appropriée quand on
regarde une orbite donnée du groupe modulaire $\PGL_2(R_v)$ sur des
irrationnels quadratiques (notons que contrairement au cas rationnel,
il y a une infinité de telles orbites). Notons $\cdot$ l'action par
homographies de $\GL_2(K_v)$ sur $\PP^1(K_v)= K_v\cup\{\infty\}$.

\medskip
\begin{theoreme}\label{theo:quadratirrat}
  Pour tout sous-groupe d'indice fini $G$ de $\GL_2(R_v)$ et tout
  irrationnel quadratique $\alpha_0\in K_v$ sur $K$, quand
  $s\to+\infty$, nous avons
$$
\frac{(q_v+1)^2\;\zeta_K(-1)\; m_0\;[\GL_2(R_v):G]}
{2\;q_v^2\;(q-1)\;|v(\operatorname{tr} g_0)|}
\;s^{-1} \sum_{\alpha\in G\cdot \alpha_0,\; h(\alpha)\leq s}
\Delta_{\alpha}\;\;\weakstar\;\; \haar_{K_v}\;,
$$
où $g_0\in G$ fixe $\alpha_0$ avec $v(\operatorname{tr} g_0)\neq 0$ et
$m_0$ est l'indice de ${g_0}^\ZZ$ dans $G_{\alpha_0}$.
\end{theoreme}

\medskip Un autre résultat d'équidistribution d'orbites d'irrationnels
quadratiques s'obtient en utilisant une complexité construite à
partir de birapports d'irrationnels quadratiques.  Nous noterons
$[a,b,c,d]=\frac{(c-a)(d-b)}{(c-b)(d-a)}$ le birapport de quatre
éléments de $K_v$ deux à deux disjoints. Si $\alpha,\beta\in K_v$ sont
deux irrationnels quadratiques sur $K$, tels que
$\alpha\notin\{\beta,\beta^\sigma\}$, notons
$$
h_\beta(\alpha)=\max\{|[\alpha,\beta,\beta^\sigma,\alpha^\sigma]|_v,\;
|[\alpha^\sigma,\beta,\beta^\sigma,\alpha]|_v\}\;,
$$
qui, comme nous allons le voir, est une autre complexité appropriée
quand $\alpha$ varie dans une orbite donnée du groupe modulaire
$\PGL_2(R_v)$ sur des irrationnels quadratiques. Par invariance de la
complexité $h_\beta$ par $\PGL_2(R_v)_\beta$, nous obtenons alors une
équidistribution vers une mesure absolument continue par rapport à la
mesure de Haar, invariante par $\PGL_2(R_v)_\beta$.

\medskip
\begin{theoreme}\label{theo:relheight}  
  Pour tout sous-groupe d'indice fini $G$ de $\GL_2(R_v)$ et tous les
  irrationnels quadratiques $\alpha_0,\beta \in K_v$ sur $K$, quand
  $s\to+\infty$, pour la convergence vague des mesures sur
  $K_v-\{\beta,\beta^\sigma\}$,
$$
\frac{(q_v+1)^2\;\zeta_K(-1)\; m_0\;[\GL_2(R_v):G]}
{2\;q_v^2\;(q-1)\;|\beta-\beta^\sigma|_v|v(\operatorname{tr} g_0)|}
\;s^{-1} \sum_{\alpha\in G\cdot \alpha_0,\; h_{\beta}(\alpha)\leq s}
\Delta_{\alpha}\;\;\weakstar\;\; 
\frac{d\haar_{K_v}(z)}{|z-\beta|_v\,|z-\beta^\sigma|_v}\;.
$$
\end{theoreme}

\medskip Le dernier résultat arithmétique affirme l'équidistribution
projective de représentations intégrales de formes normes
quadratiques vers la même mesure que dans le théorème précédent.

\medskip
\begin{theoreme}\label{theo:normform}  
  Pour tout idéal non nul $I$ de $R_v$ (de norme $N(I)$) et tout
  irrationnel quadratique $\beta \in K_v$ sur $K$, quand
  $s\to+\infty$, pour la convergence vague des mesures sur
  $K_v-\{\beta,\beta^\sigma\}$, nous avons
$$
\frac{(q_v^2-1)\;(q_v+1)\;\zeta_K(-1)\;N(I)\prod_{\ppp|I}(1+\frac{1}{N(\ppp)})}
{q_v^3\;(q-1)^2\;q^{g-1}}
\;s^{-1} \sum_{\substack{(x,y)\in R_v\times I,\;xR_v+yR_v=R_v\\|{\tt n}(x-y\beta)|_v\leq s}}
\Delta_{\frac xy}\;\;\weakstar\;\; 
\frac{d\haar_{K_v}(z)}{|z-\beta|_v\,|z-\beta^\sigma|_v}\;.
$$
\end{theoreme}

\medskip Les quatre résultats d'équidistribution ci-dessus admettent
un terme d'erreur en $\bigO(s^{-\kappa})$ pour un $\kappa>0$ quand ils
sont évalués sur une fonction localement constante à support
compact. Nous renvoyons à \cite{BroParPau16} pour des analogues aux
théorèmes \ref{theo:quadratirrat} et \ref{theo:relheight} dans
$\QQ_p$.

\section{\'Equidistribution de perpendiculaires communes
dans des arbres}
\label{sec:arbre}

Nous donnons dans cette partie l'outil géométrique principal utilisé
pour montrer les résultats de la partie précédente. Il relève de la
théorie ergodique des flots géodésiques dans les arbres. Nous
renvoyons à \cite{Serre83} pour toute information sur les actions de
groupes sur les arbres.

Soient $q\in\NN$ un entier au moins $2$, $\XX$ un arbre
$(q+1)$-régulier, d'ensemble des sommets $V\XX$ et de réalisation
géométrique $X$, et $\Aut(\XX)$ le groupe localement compact des
automorphismes sans inversion de $\XX$. Soit $\Ga$ un réseau de
$\Aut(\XX)$, c'est-à-dire un sous-groupe discret tel que
$\Ga\bs\Aut(\XX)$ admette une mesure de probabilité invariante par
$\Aut(\XX)$. Supposons que $\Ga$ préserve l'ensemble des sommets de
$\XX$ à distance paire d'un sommet donné. Pour toute partie $E$ de
$\XX$, nous noterons $\Ga_E$ le stabilisateur de $E$ dans $\Ga$.

Nous noterons $\Ga\dbs\XX$ le graphe de groupes quotient de $\XX$ par
$\Ga$ et $\vol_{\Ga\dbs\XX}$ la mesure (de masse totale notée
$\|\vol_{\Ga\dbs\XX}\|$) sur l'ensemble (discret) des sommets du
graphe quotient $\Ga\bs\XX$ définie par $\vol_{\Ga\dbs\XX}
=\sum_{[x]\in \Ga\bs V\XX} \;\frac{1}{|\Ga_x|}\;\Delta_{[x]}$.

L'espace dans lequel aura lieu l'équidistribution est l'espace (de
Bartels-Lück) localement compact $\gengeod \XX$ des géodésiques
généralisées de $\XX$, c'est-à-dire des applications continues
$\ell:\RR\ra X$ telles que $\ell(0)\in V\XX$, isométriques sur un
intervalle fermé de $\RR$, et constante sur chaque composante du
complémentaire. Il est muni de l'action de $\Aut(\XX)$ par
postcomposition et de l'action du flot géodésique (à temps $t\in\ZZ$
discret) $(t,\ell)\mapsto \{s\mapsto\ell(s+t)\}$. Il contient le
sous-espace fermé $\G\XX$ des géodésiques complètes (isométriques sur
$\RR$).

Soient $\DD^{\pm}$ deux sous-arbres de $\XX$, propres et non vides,
tels que les familles $(\ga\DD^\pm)_{\ga\in\Ga/\Ga_{\DD^\pm}}$ soient
localement finies dans $\XX$. Notons $\normalmp \DD^\pm$ le
sous-espace fermé de $\gengeod \XX$ des rayons géodésiques $\rho$
rentrant/sortant de $\DD^\pm$, c'est-à-dire isométriques sur
exactement $\mp [0,+\infty[$, avec $\rho(0)\in \DD^\pm$ et
$\rho(t)\notin\DD^\pm$ si $\mp t> 0$. Il porte une mesure borélienne
$\wt\sigma^\mp_{\DD^\pm}$ naturelle,
%(invariante par les éléments de $\Aut(\XX)$ préservant $\DD^\pm$)
dont la restriction au sous-espace des rayons géodésiques
rentrant/sortant en un point donné de $V\DD^\pm$, si non vide, est
l'unique mesure de probabilité invariante par tous les éléments de
$\Aut(\XX)$ préservant ce sous-espace.

Le résultat principal est un théorème d'équidistribution simultanée
des segments perpendiculaires communs entre $\DD^-$ et $\ga\DD^+$
lorsque $\ga$ varie dans $\Ga$. Si $\DD^-$ et $\ga\DD^+$ sont
disjoints, nous noterons $\lambda_\ga= d(\DD^-,\ga\DD^+)$ la longueur
du segment perpendiculaire commun et $\alpha^-_\ga:[0,\lambda_\ga]\ra
X$, $\alpha^+_\ga:[-\lambda_\ga,0]\ra X$ les deux paramétrages avec
$\alpha^-_\ga(0)\in \DD^-$ et $\alpha^+_\ga(0)\in \DD^+$ du segment
perpendiculaire commun, considérés comme des géodésiques généralisées
par prolongement localement constant hors de $]0,\lambda_\ga[$ et
$]-\lambda_\ga,0[\,$.

\medskip
\begin{theoreme}\label{theo:arbre} 
  Quand $t\ra+\infty$, pour la convergence vague des mesures sur
  $\gengeod \XX\times\gengeod \XX$, nous avons
$$
\frac{(q^2-1)\;(q+1)}{2\;q^2}\;\|\vol_{\Ga\dbs\XX}\|\;q^{-t}
\sum_{\substack{\ga\in\Ga/\Ga_{\DD^+}\\0<\lambda_\ga\leq t}}\;
\Delta_{\alpha^-_\ga}\otimes \Delta_{\ga^{-1}\alpha^-_\ga}
\;\;\weakstar\;\;\wt\sigma^+_{\DD^-}\otimes \wt\sigma^-_{\DD^+}\;.
$$
\end{theoreme}

Si $\Ga$ est géométriquement fini, alors ce résultat
d'équidistribution admet un terme d'erreur en $\bigO(q^{-\kappa})$
pour un $\kappa>0$ quand il est évalué sur une fonction localement
constante à support compact. D'après \cite{Paulin04b}, le groupe $\Ga$
est géométriquement fini si et seulement si le graphe de groupes
$\Ga\dbs\XX$ est réunion d'un graphe fini de groupes et d'un nombre
fini de rayons de groupes cuspidaux (c'est-à-dire dont les
stabilisateurs des arêtes orientées vers le bout fixe leur origine).

\medskip La démonstration, pour laquelle nous renvoyons à
\cite{BroParPau16}, utilise la propriété de mélange pour le carré du
flot géodésique sur le quotient par $\Ga$ de l'espace des géodésiques
complètes d'origine à distance paire d'un point base, muni de la
restriction de la mesure de Bowen-Margulis (ou mesure d'entropie
maximale). Nous montrons d'ailleurs que l'image de cette mesure sur
$\Ga\bs\G\XX$ par l'application origine $\ell\mapsto \ell(0)$ est un
multiple de la mesure $\vol_{\Ga\dbs\XX}$ sur $\Ga\bs V\XX$. Le terme
d'erreur lorsque $\Ga$ est géométriquement fini utilise la propriété
de décroissance exponentielle des corrélations du mélange, annonc\?ee 
dans \cite{Kwon15}.

\medskip Nous concluons cette note en donnant des indications sur la
déduction du théorème \ref{theo:mertens} à partir du théorème
\ref{theo:arbre}, quand, pour simplifier, $K$ est le corps des
fractions rationnelles $\FF_q(Y)$ en une variable sur $\FF_q$, la
valuation $v$ est $v_\infty:\frac{P}{Q}\mapsto \deg Q-\deg P$, de
sorte que $g=0$, $q_v=q$, $R_v=\FF_q[Y]$ et $G=\GL_2(R_v)$.

Prenons pour $\XX$ l'arbre de Bruhat-Tits de $(\PGL_2,K_v)$ (voir
\cite{Serre83}), dont les sommets sont les classes d'homothéties de
$\OOO_v$-réseaux de $K_v\times K_v$, dont l'espace des bouts
$\partial_\infty\XX$ s'identifie avec $\PP_1(K_v)=K_v\cup\{\infty\}$,
de sorte que l'ensemble des extrémités en $+\infty$ des géodésiques
complètes dont l'extrémité en $-\infty$ est le point $\infty\in
\PP_1(K_v)$ et qui passent par $*_v=[\OOO_v\times\OOO_v]$ soit
exactement $\OOO_v$.

Prenons pour $\Ga$ le réseau de Nagao $\PGL_2(R_v)$ (voir par exemple
\cite{Weil70}), dont le graphe de groupes quotient $\Ga\dbs\XX$ est
formé d'un rayon cuspidal recollé en son origine à une arête de
groupes, et dont le covolume $\|\vol_{\Ga\dbs\XX}\| =
\frac{2}{(q-1)(q^2-1)}$ est bien connu.
%dont le groupe du sommet d'origine est
%$\PGL_2(\FF_q)$, et dont le groupe du sommet à distance $n\neq 0$ de
%l'origine est d'ordre $(q-1)q^{n}$, de sorte que
%$\vol_{\Ga\dbs\XX}=\frac{2}{(q-1)(q^2-1)}$.

Prenons pour sous-arbres $\DD^-$ et $\DD^+$ tous deux l'horoboule
$\H_\infty$ de $\XX$ centrée en $\infty\in \PP_1(K_v)$, dont le bord
passe par $*_v$, de sorte que si $\beta_\infty(x,y)= \lim_{z\ra
  \infty} d(x,z)-d(y,z)$, alors $V\DD^-=V\DD^+=\{x\in V\XX\;:\;
\beta_\infty(x,*_v)\leq 0\}$.  Notons que si $\theta$ est
l'application continue et propre qui à un rayon géodésique de
$\normalout \DD_-$ associe son point l'infini, alors
$\theta_*(\wt\sigma^+_{\DD^-})= \haar_{K_v}$.

Un petit calcul montre que l'image par $\ga$ de $\DD^+$ est
l'horoboule de $\XX$ centrée en $\ga\infty=\frac{a}{c}$ avec $a,b\in
R_v$ premiers entre eux, et dont le segment de perpendiculaire commun
avec $\DD^-$ est de longueur $-2\,v(c)=2\frac{\ln |c|_v}{\ln
  q_v}$. L'application prolongeant ce segment en un rayon géodésique
de point à l'infini $\ga\infty$ étant uniformément continue, par le
changement de variable $s=q_v^{\;\frac{t}{2}}$ et par la continuité de
$\theta_*$ pour les topologies vagues, le théorème \ref{theo:mertens}
découle du théorème \ref{theo:arbre} (en ne considérant que la
première composante).

Les théorèmes \ref{theo:quadratirrat}, \ref{theo:relheight} et
\ref{theo:normform} s'obtiennent de même en prenant $\DD^-=\H_\infty$
et $\DD^+=\; ]\alpha_0,\alpha_0^\sigma[$ la géodésique de $\XX$ de
points à l'infini $\alpha_0,\alpha_0^\sigma$ pour le premier,
$\DD^-=\; ]\beta,\beta^\sigma[$ et $\DD^+=\; ]\alpha_0,
\alpha_0^\sigma[$ pour le second, et $\DD^-=\; ]\beta,\beta^\sigma[$
et $\DD^+=\; \H_\infty$ pour le troisième.

{\small\

\noindent \begin{tabular}{l}
Laboratoire de math\'ematique d'Orsay, UMR 8628 Universit\'e Paris-Sud et CNRS\\
Universit\'e Paris-Saclay, 91405 ORSAY Cedex, FRANCE\\
{\it e-mail: anne.broise@math.u-psud.fr}
\end{tabular}

\medskip

\noindent \begin{tabular}{l} 

Department of Mathematics and Statistics, P.O. Box 35\\ 
40014 University of Jyv\"askyl\"a, FINLAND.\\
{\it e-mail: jouni.t.parkkonen@jyu.fi}
\end{tabular}

\medskip
\noindent \begin{tabular}{l} Laboratoire de math\'ematique d'Orsay,
  UMR 8628 Universit\'e Paris-Sud et CNRS\\ Universit\'e Paris-Saclay,
  91405 ORSAY Cedex, FRANCE\\ {\it e-mail:
    frederic.paulin@math.u-psud.fr}
\end{tabular}

}

\end{document}